\documentclass[10pt,a4paper]{article}
\usepackage[francais,english]{babel}
\usepackage{amsmath}
\usepackage{amssymb}

\newcommand{\HH}{\ensuremath{\mathcal{H}}}
\newcommand{\mcc}{\mathbb{C}}
\newcommand{\rHH}{{\rm HH}}
\newcommand{\Tor}{{\rm Tor}}
\newcommand{\ib}{\item[$\bullet$]}
\newcommand{\bT}{\mathbf{T}}
\newcommand{\bb}{\mathbf{b}}
\newcommand{\ba}{\mathbf{a}}
\newcommand{\hh}{\mathbf{h}}
\newcommand{\ehr}[1]{\spadesuit_{q,t}(#1)}
\newcommand{\dimgr}{{\rm dimgr}}
\newcommand{\rH}{{\rm H}}
\newcommand{\ty}{{$\newline$ }}
\newcommand{\lrv}[1]{\langle#1\rangle_V}

\newtheorem{theo}{Theorem}
\newtheorem{propo}{Proposition}

\newtheorem{lemm}{Lemma}
\newtheorem{remm}{Remark}
\newtheorem{fact}{Fact}

\begin{document}
\title{Markov property and Khovanov-Rozansky homology: Coxeter Case.}
\author{Trafim LASY}
\date{\today}
\maketitle
\begin{abstract}
In \cite{Kh07} Khovanov gave a construction of Khovanov-Rozansky homology, new triple-graded link invariant, taking the Hochschild homology
of the terms in Rouquier complexes. Its bigraded Euler characteristic provides the Markov trace of type $A_n$,
which, in its turn, induces the well-known HOMFLY polynomial. Gomi \cite[(3.3)]{Go06} has extended  the definition of Markov trace
on a Hecke algebra to all Coxeter groups. One can see that Khovanov's and Rouquier's constructions can be directly extended to these
 groups. We give a detailed proof that the corresponding Euler characteristic provides a Markov trace in the sense
of Gomi for any Coxeter group.
\end{abstract}
\bigskip
Let $(W,S)$ be a Coxeter system with $|S|=n<\infty$. Let $V$ be the $\mcc$-vector space with the basis $\{e_s\}_{s\in S}$ carrying
the geometric representation of $W$. Then we have a natural action of $W$ on the symmetric algebra $R=S(V^*)$, which will be
regarded as the graded polynomial ring $\mcc[\{x_s\}_{s\in S}]$ with $\deg x_s=1$.

Now we briefly recall the notions of Rouquier complexes needed to give Khovanov's construction of Khovanov-Rozansky
homology. Let $B_W=\left<\sigma_s\right>_{s\in S}$ be the corresponding to $W$ braid group. For $s\in S$ define the graded
$R$-bimodule $B_s = R\otimes_{R^s}R$, where $R^s$ is the ring of polynomials invariant under $s$. Graduation on $B_s$ is induced by the
graduation on $R$. To the braid generator $\sigma_s$ (resp. $\sigma_s^{-1}$) Rouquier associates cochain complexes of graded
$R$-bimodules:
$$
F(\sigma_s): 0\rightarrow B_s \stackrel{m}{\rightarrow} R \rightarrow 0 \quad\mbox{resp.}\quad
F(\sigma_s^{-1}) : 0\rightarrow R \stackrel{\eta}{\rightarrow} B_s(1) \rightarrow 0
$$
where $B_s(1)$ means the shift in the graduation by $1$, $m(a\otimes b) = ab$ and $\eta(a) = a\otimes x_s + ax_s\otimes 1$. Bimodules
$B_s$ and $B_s(1)$ are placed in cohomological degree $0$. Rouquier \cite[9.]{Ro06} proves the following
\begin{theo}[Rouquier]
Map $\sigma_s^{\pm 1}\to F(\sigma_s^{\pm 1})\otimes_R -$ induces
a well-defined action of the braid group $B_W$ on the category
of cochain complexes of graded $R$-bimodules up to homotopy.
\end{theo}
The above theorem allows to define the {\em Rouquier complex} $F(\sigma)$ of $\sigma\in B_W$ as the image
of the complex $0\rightarrow R\rightarrow 0$ by the action of $\sigma$.

Recall that the {\em Hochschild homology} $\rHH(R,M)$ of an $R$-bimodule $M$ is defined as the direct sum of $\Tor_i^{R\otimes R}(R,M)$.
For $\sigma\in B_W$ Khovanov takes the Hochschild homology of each term in the Rouquier complex $F(\sigma)$
and gets a complex of bigraded spaces
\begin{equation*}
\ldots \xrightarrow[]{\rHH(\partial)} \rHH(R,F^j(\sigma))\xrightarrow[]{\rHH(\partial)}
 \rHH(R,F^{j+1}(\sigma))\xrightarrow[]{\rHH(\partial)}\ldots
\end{equation*}
The homology of the above complex is called {\em Khovanov-Rozansky homology}. When $W$ is the symmetric group Khovanov
proves the following
\begin{theo}[Khovanov]
The triple-graded dimension of Khovanov-Rozansky homology is an invariant of oriented links. Its bigraded
Euler characteristic provides the HOMFLY $2$-variable polynomial.
\end{theo}
Let $\HH_W:=\mcc(q)B_W/\left<(\sigma_s-q)(\sigma_s+1)\right>$ be the Hecke algebra of $W$.
We will explain how the Euler characteristic of Khovanov-Rozansky homology provides a trace function on $\HH_W$.
For $\sigma\in B_W$ define an element $\lrv\sigma\in\mcc(q,t)$ as
\begin{equation*}
\left<\sigma\right>_{V} =\left(\frac{1-q}{1+tq}\right)^n\cdot \sum_{i,j,k}(-1)^i
\dim_{\mcc}\rHH_j(R,F^i(\sigma))_k\cdot q^k\cdot t^j,
\end{equation*}
where the index $i$ corresponds to the cohomological grading of the Rouquier complex,
$j$ - to the Hochschild homology grading and $k$ - to the polynomial
grading of $R$-modules $\rHH_j(R,F^i(\sigma))$. By linearity we get a $\mcc(q)$-linear map
$\lrv\cdot:\mcc(q)B_W\rightarrow\mcc(q,t)$. The following proposition is well-known (see e.g. \cite[(4.4.6)]{La12}):
\begin{propo} \label{01}
The map $\lrv\cdot$ has the following properties:
\begin{itemize}
\ib $\left<1\right>_V=1$ (Normalizing property)
\ib $\left<ab\right>_V=\left<ba\right>_V$ for any $a, b\in \mcc(q)B_W$ (Trace property)
\ib  For any $a\in\mcc(q)B_W$ and $s\in S$ we have $\left<a\sigma_s\right>_V = q\left<a\sigma_s^{-1}\right>_V +
(q-1)\left<a\right>_V$ (Hecke property)
\end{itemize}
\end{propo}
This proposition implies that the map $\lrv\cdot$ can be factorized through the Hecke algebra $\HH_W$. It induces a
trace $\tau_{kr}$ on $\HH_W$ called {\em Khovanov-Rozansky trace}.

Let $\bT_s$ be the image of $\sigma_s$ in $\HH_W$. Recall, that Gomi defines a {\em Markov trace} with a parameter $z\in\mcc(q,t)$ on
$\HH_W$ as a $\mcc(q)$-linear function $\tau$ with the properties
\begin{itemize}
\ib $\tau(1)=1$, $\tau(\ba\bb)=\tau(\bb\ba)$ for any $\ba, \bb\in\HH_W$
\ib $\tau(\hh\bT_s)=z\tau(\hh)$ for any $s\in S$ and any $\hh$ from the parabolic subalgebra $\HH_{W, S\backslash\{s\}}$
of $\HH_W$ generated by the set $\{\bT_{s'}\}_{s'\in S\backslash\{s\}}$. This property is called {\em Markov property}.
\end{itemize}
The result of this article is the following theorem:
\begin{theo} \label{02}
The  trace $\tau_{kr}$ is a Markov trace with the parameter $z=\frac{tq(q-1)}{tq+1}$.
\end{theo}
{\bf Proof.} The idea of the proof is based on the proof of \cite[Th. 1]{Kh07} where Koszul complexes are used.
Let $e_1,\cdots,e_n$ be any basis of the vector space $V$. We denote the dual $e^*_i$ by $x_i\in R$. We will use the following notations:
$$
\begin{array}{llll}
R_x^n := \mcc[x_1,\ldots,x_n], & R_y^n := \mcc[y_1,\ldots,y_n], & R_z^n := \mcc[z_1,\ldots,z_n], &
R_{x,y,z}^n := R_x^n\otimes_\mcc R_y^n\otimes_\mcc R_z^n, \\
R_{x,y}^n := R_x^n\otimes_\mcc R_y^n, & R_{x,z}^n := R_x^n\otimes_\mcc R_z^n, & R_{y,z}^n := R_y^n\otimes_\mcc R_z^n.
\end{array}
$$
The action of $W$ on $R_y^n,R_z^n$ is exactly the same as on $R = R_x^n$ and is directly induced on
$R_{x,y}^n,R_{x,z}^n,R_{y,z}^n,R_{x,y,z}^n$.

Let $A$ be a graded ring and $a_1,\ldots a_k$ some homogeneous elements in $A$. We denote by $(a_1,\ldots,a_k)_A$
the corresponding {\em Koszul complex} which is defined as the tensor product over $A$ of small complexes
$$
0 \rightarrow A(-\deg a_i)\stackrel{\times a_i}{\longrightarrow} A \rightarrow 0, \quad
\mbox{where $\times a_i$ means the multiplication by $a_i$.}
$$
\begin{fact}
Let $(a_1,\ldots,a_k)_A$ be a Koszul complex and $\lambda\in A$ be an element such that $\lambda a_i$ is homogeneous
of the same degree as $a_j$ for some $i\neq j$. Then we have the following isomorphism of complexes:
$$(a_1,\ldots,a_i,\ldots,a_j\ldots,a_k)_A\simeq(a_1,\ldots,a_i,\ldots,a_j+\lambda a_i\ldots,a_k)_A.$$
\end{fact}
\begin{fact}
Let $A,B$ be some $\mcc$-algebras, $K^\bullet$ some complex of $A\otimes_\mcc B$-modules
and $L^\bullet$ a complex of $A$-modules. Then there is the obvious isomorphism of complexes
$K^\bullet\otimes_A L^\bullet \simeq K^\bullet\otimes_{A\otimes B} (L^\bullet\otimes_\mcc B)$.
\end{fact}
Proposition \ref{01}, the definition and the $\mcc(q)$-linearity of $\tau_{kr}$ imply that in order to prove the theorem
it is enough to show that $\lrv{b\sigma_s} = z\lrv b$
for any $s\in S$ and $b$ being a product of some generators $\sigma_t, t\neq s$ of the braid group $B_W$.

For a bigraded space $M$ with finite-dimensional homogeneous components we denote its bigraded dimension by $\dimgr_{q,t}M$.
For a finitely generated graded $R$-bimodule $N$ we denote by $\ehr{N}$ the bigraded dimension of its Hochschild homology
$\dimgr_{q,t}\rHH(R,N)=\sum_{j,k}\dim_\mcc\rHH_j(R,N)_k\cdot q^k\cdot t^j$. Recall that by definition
\begin{align*}
&\left(\frac{1+tq}{1-q}\right)^n\lrv{b\sigma_s} = \sum_i(-1)^i\ehr{F^i(b\sigma_s)} =
 \sum_i(-1)^i\ehr{F^i(b)\otimes_R B_s} - \sum_i(-1)^i\ehr{F^i(b)} \\
&\left(\frac{1+tq}{1-q}\right)^nz\lrv b = \frac{tq(q-1)}{tq+1}\sum_i(-1)^i\ehr{F^i(b)}.
\end{align*}
Thus it is enough to prove the following equality: $(tq+1)\ehr{\Theta_b\otimes_R B_s} = (tq^2+1)\ehr{\Theta_b}$, \ty
where the bimodule $\Theta_b$ is a tensor product over $R$ of some $B_t$'s ($t\neq s$) since every $F^i(b)$ is a direct sum of
 such bimodules.

Let us look at our vector space $V$. Since $\dim_\mcc V = n > n-1 = |S\setminus\{s\}|$ there exists a vector $v\in V^*$ invariant under
the action of any $s'\in S, s'\neq s$. We can assume that $v = x_n$ and $x_1,\ldots,x_{n-1}\in V^s$.
Let $u_x = x_n - w_x$ $(w_x\in V^{*s})$
 be an element of $V^*$ dual to the root of $s$ (vectors $u_y = y_n - w_y$ and $u_z = z_n - w_z$ will be the corresponding copies of
$u_x$ in $R_y^n$ and $R_z^n$).
\begin{lemm}
There is an isomorphism of $R_{x,y}^n$-modules $\Theta_b\simeq\Theta_b'\otimes_\mcc\mcc[x_n]$,
where $\Theta_b'$ is some $R_{x,y}^{n-1}$-module, $y_n$ acts on $\mcc[x_n]$ by multiplication by $x_n$.
\end{lemm}
{\bf Proof.} Easy consequence from the fact that $\Theta_b$ is a tensor product of some $B_t$'s. $\Box$
\begin{propo}
Let $\widetilde{\Theta_b'}$ be a free $R_{x,y}^{n-1}$-resolution of $\Theta_b'$. Then the following complex $K_{b,s}$ is
a free $R_{x,y,z}^n$-resolution of $\Theta_b\otimes_{R_y^n} B_s$:
$$
\left(\widetilde{\Theta_b'}\otimes_\mcc\mcc[x_n,y_n]\otimes_\mcc R_z^n\right) \otimes_{R_{x,y,z}^n}
(x_n-y_n,y_1-z_1,\ldots,y_{n-1}-z_{n-1},u_y^2-u_z^2)_{R_{x,y,z}^n}.
$$
\end{propo}
{\bf Proof.} Indeed, the above complex is isomorphic to
\begin{align*}
&\left[\left(\widetilde{\Theta_b'}\otimes_\mcc\mcc[x_n,y_n]\right) \otimes_{R_{x,y}^n}
(x_n-y_n)_{R_{x,y}^n}\right]\otimes_{R_{x,y}^n} (y_1-z_1,\ldots,y_{n-1}-z_{n-1},u_y^2-u_z^2)_{R_{x,y,z}^n} \\
\simeq &\left[\widetilde{\Theta_b'}\otimes_{R_{x,y}^{n-1}}(x_n-y_n)_{R_{x,y}^n}\right]
\otimes_{R_y^n} (y_1-z_1,\ldots,y_{n-1}-z_{n-1},u_y^2-u_z^2)_{R_{y,z}^n}
\end{align*}
Since the complex $\widetilde{\Theta_b'}\otimes_{R_{x,y}^{n-1}}(x_n-y_n)_{R_{x,y}^n}
$is isomorphic to $\widetilde{\Theta_b'}\otimes_\mcc(x_n-y_n)_{\mcc[x_n,y_n]}$, it is a free $R_{x,y}^n$-resolution of $\Theta_b$.
At the same time the complex $(y_1-z_1,\ldots,y_{n-1}-z_{n-1},u_y^2-u_z^2)_{R_{y,z}^n}$ considered as
a complex of $R_y^n$-modules is homotopy equivalent to the free $R_y^n$-module $B_s$. This proves that their tensor product over
$R_y^n$ is a free $R_{x,y,z}^n$-resolution of $\Theta_b\otimes_{R_y^n} B_s$. $\Box$ \ty
In order to calculate $\rHH(R,\Theta_b\otimes_R B_s)$ we take the tensor product over $R_{x,z}^n$ of
$K_{b,s}$ and the free $R_{x,z}^n$-resolution of $R$ represented by the Koszul complex
$(x_1-z_1,\ldots,x_n-z_n)_{R_{x,z}^n}$. This tensor product equals
$$
K_{b,s}\otimes_{R_{x,y,z}^n}(x_1-z_1,\ldots,x_n-z_n)_{R_{x,y,z}^n},\quad \mbox{which, in its turn, is equal to }
$$
\begin{align*}
&\left(\widetilde{\Theta_b'}\otimes_\mcc\mcc[x_n,y_n]\otimes_\mcc R_z^n\right)
\otimes_{R_{x,y,z}^n}(x_n-y_n,y_1-z_1,\ldots,y_{n-1}-z_{n-1},u_y^2-u_z^2,x_1-z_1,\ldots,x_n-z_n)_{R_{x,y,z}^n} \\
&\simeq\left(\widetilde{\Theta_b'}\otimes_\mcc\mcc[x_n,y_n]\otimes_\mcc R_z^n\right)
\otimes_{R_{x,y,z}^n}(x_n-y_n,y_1-z_1,\ldots,y_{n-1}-z_{n-1},0,x_1-z_1,\ldots,x_n-z_n)_{R_{x,y,z}^n},
\end{align*}
where the element $u_y^2-u_z^2 = (y_n-z_n-w_y+w_z)(y_n+z_n-w_y-w_z)$ of homogeneous degree $2$ was "killed" by elements
$x_n-y_n,x_n-z_n,y_1-z_1,\ldots,y_{n-1}-z_{n-1}$. Denoting the homology of the complex
$$
\left(\widetilde{\Theta_b'}\otimes_\mcc\mcc[x_n,y_n]\otimes_\mcc R_z^n\right)
\otimes_{R_{x,y,z}^n}(x_n-y_n,y_1-z_1,\ldots,y_{n-1}-z_{n-1},x_1-z_1,\ldots,x_n-z_n)_{R_{x,y,z}^n}
$$
by $\rH_{b,s}^\bullet$ (which is double-graded) we have $\ehr{\Theta_b\otimes_R B_s} = (1+tq^2)\dimgr_{q,t}\rH_{b,s}^\bullet$.

The proof that $\ehr{\Theta_b} = (1+tq)\dimgr_{q,t}\rH_{b,s}^\bullet$
is completely analogous. First we make a formal note that the complex $L_{b,s}$:
$$\left(\widetilde{\Theta_b'}\otimes_\mcc\mcc[x_n,y_n]\otimes_\mcc R_z^n\right) \otimes_{R_{x,y,z}^n}
(x_n-y_n,y_1-z_1,\ldots,y_{n-1}-z_{n-1},y_n-z_n)_{R_{x,y,z}^n}$$
is a free $R_{x,y,z}^n$-resolution of $\Theta_b$. Then we take the tensor product over $R_{x,z}^n$ of
$L_{b,s}$ and the Koszul complex $(x_1-z_1,\ldots,x_n-z_n)_{R_{x,z}^n}$.
And we get
\begin{align*}
&\left(\widetilde{\Theta_b'}\otimes_\mcc\mcc[x_n,y_n]\otimes_\mcc R_z^n\right)
\otimes_{R_{x,y,z}^n}(x_n-y_n,y_1-z_1,\ldots,y_{n-1}-z_{n-1},y_n-z_n,x_1-z_1,\ldots,x_n-z_n)_{R_{x,y,z}^n} \\
&\simeq\left(\widetilde{\Theta_b'}\otimes_\mcc\mcc[x_n,y_n]\otimes_\mcc R_z^n\right)
\otimes_{R_{x,y,z}^n}(x_n-y_n,y_1-z_1,\ldots,y_{n-1}-z_{n-1},0,x_1-z_1,\ldots,x_n-z_n)_{R_{x,y,z}^n},
\end{align*}
where the element $y_n - z_n$ of homogeneous degree $1$ was "killed" by elements
$x_n-y_n$ and $x_n-z_n$. This ends the proof. $\Box$
\begin{remm}
Rapha\"el Rouquier told me in private communication that he has a "categorified" version of theorem \ref{02}
which will be published later.
\end{remm}

\end{document}